\documentclass[12pt]{article}
\usepackage{mathtools}
\usepackage[margin=1in]{geometry}
\usepackage{amsfonts}
\usepackage{amssymb}
\usepackage{tikz}
\usepackage{graphicx}
\usepackage{fancyhdr}
\usetikzlibrary{arrows}
\usetikzlibrary{shapes}

         \def \s {\sigma}     \def \tss {\textsuperscript}     \def \sse {\subseteq}      \def \rar {\rightarrow}               
 \def \td {\text{-}} \def \fbe {f_\text{before}} \def \faf {f_\text{after}} \def \fen {f_\text{end}} \def \fbu {f_\text{bump}}
\begin{document}
\begin{center}
A New Quantity Counted by OEIS Sequence A006012
\\
Yonah Biers-Ariel
\\
Rutgers University
\end{center}
\noindent We prove a conjecture of Callan in \cite{oeis} that OEIS sequence A006012 counts a certain kind of permutation. Call this sequence $(a_n)_{n=1}^\infty$; then $a_n$ is defined by $a_1=1$, $a_2=2$, and $a_n=4a_{n-1}-2a_{n-2}$ (the actual sequence in the OEIS is offset by one, so $a_0=1$, $a_1=2$, and the recursion is the same). The conjecture states that $a_n$ is equal to the number of permutations of length $n$ for which no subsequence $abcd$ has the following two properties: $c= b+1$ and $\max\{a,c\} < \min\{b,d\}$. 

We can rewrite this conjecture in the language of pattern avoidance, in particular, using the dashed notation for generalized pattern avoidance introduced in \cite{Babson}. Therefore, we define a pattern to be a permutation $\pi_1...\pi_k$, some (or all) of whose elements may be separated by dashes. We say that a subsequence of a permutation is an occurrence of a pattern if (i) all the elements have the same relative order as the elements of the pattern, and (ii) if there is no dash between the $i\tss{th}$ and $i+1\tss{th}$ elements of the pattern, then the $i\tss{th}$ and $i+1\tss{th}$ element of the subsequence occur consecutively in the permutation. We say that a permutation avoids a pattern if it does not contain any occurrence of the pattern, and a permutation avoids a set of patterns if it does not contain any occurrence of any of them. If $A$ is a set of patterns, we will write Av$(A)$ for the set of permutations which avoid them, and $\text{Av}_n(A)$ for the set of length-$n$ permutations which avoid them. The following two examples should help clarify these definitions.
\\
\\
\textbf{Example:} The permutation 251346 contains the subsequence 5146 which is an occurrence of the pattern 3-1-24 because the elements of the subsequence occur in the same relative order as 3124, and the 4 and 6 are consecutive in the original permutation (the 5 and 1 are also consecutive - that is allowed but not necessary). 
\\
\\
\textbf{Example:} The permutation 251346 avoids 32-1-4 (i.e. $251346 \in \text{Av}_n(\{32\td1\td4\}) \sse \text{Av}(\{32\td1\td4\})$). 
\\
\\
Using this notation, we rewrite the conjecture as $a_n = |\text{Av}_n(\{1\text{-}32\td4, 1\td42\td3, 2\td31\td4, 2\td41\td3\}|$.

We will prove two propositions. The first is that if $A= \{1\text{-}32\td4, 1\td42\td3, 2\td31\td4, 2\td41\td3\}$ and $B=\{1\td3\td2\td4, 1\td4\td2\td3, 2\td3\td1\td4, 2\td4\td1\td3\}$, then Av$(A)$ and Av$(B)$ are the same set. The second proposition is that $|\text{Av}_n(B)|$ follows the defining recurrence of $a_n$, i.e. $|\text{Av}_1(B)|=1, |\text{Av}_2(B)|=2, |\text{Av}_n(B)|=4|\text{Av}_{n-1}(B)|-2|\text{Av}_{n-2}(B)|$.
\\
\\
\noindent \textbf{Proposition 1:} Let $A= \{1\text{-}32\td4, 1\td42\td3, 2\td31\td4, 2\td41\td3\}$ and $B=\{1\td3\td2\td4, 1\td4\td2\td3, 2\td3\td1\td4, 2\td4\td1\td3\}$. The sets Av$(A)$ and Av$(B)$ are the same.
\\
\\
\textit{Proof:} 
We will show that any permutation containing an occurrence of an element of $B$ must also contain an occurrence of an element of $A$ (the converse is immediately clear). Let $\pi$ be a permutation. First, note that a subpermutation $\pi_a\pi_b\pi_c\pi_d$ of $\pi$ is an occurrence of a pattern in $A$ if and only if $c=b+1$ and $\max\{\pi_a,\pi_c\} < \min\{\pi_b,\pi_d\}$ (in fact, this is the definition Callan provides in the OEIS). Similarly, a subpermutation $\pi_a\pi_b\pi_c\pi_d$ of $\pi$ is an occurrence of a pattern in $B$ if and only if $\max\{\pi_a,\pi_c\} < \min\{\pi_b,\pi_d\}$.

Choose any element of $B$, and suppose that $\pi$ contains an occurrence of this element. As noted above, this means that we can find $a <b<c<d$ such that $\max\{\pi_a,\pi_c\} < \min\{\pi_b, \pi_d\}$. Let $e$ be the largest index less than $c$ such that $\pi_e > \max\{\pi_a, \pi_c\}$, i.e. $e= \max\{i : i <c, \pi_i > \max\{\pi_a, \pi_c\}\}$. Because $b$ is an element of $\{i : i <c, \pi_i > \max\{\pi_a, \pi_c\}\}$, it follows that $e$ exists and $a < b \le e < e+1\le c< d$. Now, we claim that $\pi_a\pi_e\pi_{e+1}\pi_d$ is an occurrence of a pattern in $A$. Obviously $e +1 =e+1$, and so it remains to check that $\max\{\pi_a,\pi_{e+1}\} < \min\{\pi_e, \pi_d\}$. Because $\max\{\pi_a,\pi_c\} < \min\{\pi_b, \pi_d\}$, we conclude that $\pi_a < \pi_d$ and by the choice of $e$, we also have $\pi_a < \pi_e$. Now, either $e+1 = c$, in which case $\pi_{e+1}=\pi_c$, or else $\pi_{e+1} < \max\{\pi_a, \pi_c\}$ because otherwise we would have chosen $e+1$ as the $\max\{i : i <c, \pi_i > \max\{\pi_a, \pi_c\}\}$ instead of $e$. It follows that $\pi_{e+1} \le \max\{\pi_a,\pi_c\} < \pi_d, \pi_e$ for the same reasons as $\pi_a$. Therefore, $\max\{\pi_a,\pi_{e+1}\} < \min\{\pi_e, \pi_d\}$ and $\pi_a\pi_e\pi_{e+1}\pi_d$ is an occurrence of a pattern in $A$. We conclude that the permutations avoiding the patterns of $A$ are the same as the permutations avoiding the patterns of $B$.
\\
\\
\noindent \textbf{Proposition 2:} The number of permutations of length $n$ avoiding all patterns in $B$ (and hence in $A$) satisfies the recurrence $a_1=1, a_2=2, a_n = 4a_{n-1} - 2a_{n-2}$. 
\\
\\
\textit{Proof:} Since $\text{Av}_{1}(B) = \{1\}$ and $\text{Av}_{2}(B) = \{12, 21\}$, the initial conditions hold. Our strategy will be as follows: given $\text{Av}_n(B)$, define four maps which, when all of them are applied to all the permutations of $\text{Av}_{n-1}(B)$, will generate all of the permutations of $\text{Av}_n(B)$. Then, we will count how many permutations of $\text{Av}_n(B)$ are double counted in this way, and find that there are two for every element of $\text{Av}_{n-2}(B)$, thereby establishing the recurrence.

Note that, for a permutation to avoid all patterns of $A$, it must be the case that either 1 and 2 occur consecutively (not necessarily in that order) or either 1 or 2 is the last element of the permutation. This observation motivates the following definitions of the four maps $\fbe, \faf, \fen, \fbu$. Let $\fbe$ be the function that inputs a permutation and outputs that permutation with all elements increased by 1 and a 1 inserted immediately before the new 2. Let $\faf$ be the function that also inputs a permutation and outputs that permutation with all the elements increased by 1 and a 1 inserted immediately after the new 2. Similarly, let $\fen$ be the function that inputs a permutation, increases all its elements by 1 and puts a 1 at the end of it, and let $\fbu$ be the function that inputs a permutation, increases all its elements by 1, replaces the new 2 with a one and puts a 2 at the end. The following example gives a concrete illustration of the four functions.
\\
\\
\textbf{Example:} Let $\pi= 31542$. Then $\fbe(\pi) = 412653, \allowbreak \faf(\pi) = 421653, \allowbreak \fen(\pi) = 426531, \allowbreak \text{ and } \fbu(\pi) = 416532$. Note that $\pi \in \text{Av}(B)$, and so are all its images.
\\
\\
We claim that (i) these four functions all map elements of $\text{Av}_{n-1}(B)$ to elements of $\text{Av}_n(B)$ and (ii) $\fbe(\text{Av}_{n-1}(B)) \cup \faf(\text{Av}_{n-1}(B)) \cup \fen(\text{Av}_{n-1}(B)) \cup \fbu(\text{Av}_{n-1}(B)) \supseteq \text{Av}_{n}(B)$ (by claim (i), we could replace the `$\supseteq$' in claim (ii) with `$=$'). To verify the first claim, choose some $\s \in \text{Av}_{n-1}(B)$, and consider each function in turn. If $\fbe(\s)$ or $\faf(\s)$ contains an occurrence of a pattern in $A$, then this occurrence must use no more than 1 of the elements 1 and 2 (because they are consecutive in both $\fbe(\s)$ and $\faf(\s)$ but can't be in any pattern in $A$). Therefore, either this occurrence fails to use 1 and would have already been an occurrence of the pattern in $\s$, or else it fails to use 2, in which case it could have used 2 instead of 1 and been an occurrence of the pattern in $\s$. Thus, no such occurrence is possible in $\fbe(\s)$ or $\faf(\s)$. In addition, if $\fen(\s)$ or $\fbu(\s)$ contains an occurrence of a pattern in $A$, then this occurrence cannot use the last element because that element is either a 1 or a 2, and patterns in $A$ only end with $3$ or $4$. So, this occurrence would already be an occurrence of the pattern in $\s$, and therefore cannot exist.

To verify the second claim, chose some $\pi \in \text{Av}_{n}(B)$. As previously noted, either 1 and 2 occur consecutively in $\pi$, or else either 1 or 2 is the final element of $\pi$. Let $\pi'$ be $\pi$ with the 1 removed and each element decreased by 1. We have introduced no new patterns, and so $\pi' \in \text{Av}_{n-1}(B)$. Suppose that 1 occurs immediately before 2 in $\pi$, then $\fbe(\pi') =\pi$. If the 1 occurs immediately after 2 in $\pi$, then $\faf(\pi')= \pi$. If the 1 occurs at the end of $\pi$, then $\fen(\pi')=\pi$. If the 2 occurs at the end of $\pi$, then we will need to define $\pi''$ which is $\pi$ with the 1 removed, the 2 moved the position where the 1 used to be, and each element decreased by 1. Again, we have introduced no new patterns, and so $\pi'' \in \text{Av}_{n-1}(B)$, and $\fbu(\pi'')=\pi$.

If these four functions all had disjoint ranges, we could conclude that $a_n= 4a_{n-1}$. Unfortunately, some permutations are counted twice. Each $f$ outputs a certain kind of permutation: $\fbe$ outputs permutations where 1 immediately precedes 2, $\faf$ outputs permutations where 2 immediately precedes 1, $\fen$ outputs permutations where 1 occurs at the end, and $\fbu$ outputs permutations where 2 occurs at the end. If a permutation fulfills two of these criteria it will be double-counted. Such permutations must be counted once by either $\fbe$ or $\faf$ and again by either $\fen$ or $\fbu$ because no permutation can be counted by both $\fbe$ and $\faf$ or both $\fen$ and $\fbu$. Thus, the final two elements of such permutations are 1 and 2 (not necessarily in that order). Let $g : \text{Av}_{n}(B) \rar \text{Av}_{n-2}(B)  $ be defined as the function which takes a permutation, removes from it the elements 1 and 2, and reduces all other elements by 2. If we restrict $g$ to those permutations which end in either 12 or 21, $g$ becomes a 2-to-1 map from the double-counted permutations of $\text{Av}_{n}(B)$ to the permutations of $\text{Av}_{n-2}(B)$, and so the number of double-counted permutations is twice $a_{n-2}$. It follows that $a_n=4a_{n-1}-2a_{n-2}$. 
\\
\\
The author would like to thank his (intended) advisor Doron Zeilberger for introducing him to the conjecture and fixing typos in the original draft.

\end{document}